\documentclass{amsart}
\usepackage{graphicx}
\theoremstyle{definition}

\theoremstyle{remark}

\numberwithin{equation}{section}

\begin{document}

\title{Optimal transportation and monotonic quantities on evolving manifolds}%
\author{Hong Huang}%
\address{School of Mathematical Sciences,Key Laboratory of Mathematics and Complex Systems,
Beijing Normal University,Beijing 100875, P. R. China}%
\email{hhuang@bnu.edu.cn}%

\thanks{Partially supported by NSFC no.10671018.}%
\subjclass{53C44}%

\keywords{optimal transportation, $\mathcal{L}$-length, Boltzmann-Shannon entropy, evolving manifolds}%

\begin{abstract}

In this note we will adapt Topping's $\mathcal{L}$-optimal
transportation theory for Ricci flow to a more general situation,
i.e. to a closed manifold $(M,g_{ij}(t))$ evolving by
$\partial_tg_{ij}=-2S_{ij}$, where $S_{ij}$ is a symmetric tensor
field of (2,0)-type on $M$. We extend  some recent results of
Topping, Lott and Brendle, generalize the monotonicity of List's
(and hence also of Perelman's) $\mathcal{W}$-entropy, and recover
the monotonicity of M$\ddot{u}$ller's (and hence also of Perelman's)
reduced volume.

\end{abstract} \maketitle

\section {Introduction}

Since Monge introduced the optimal transportation problem, many
beautiful works
 have been done, in particular in the last several decades. For an extensive discussion see Villani [V].
 Recently,
Topping, Lott, Brendle and some other authors considered this
problem on a manifold evolving according to Hamilton's Ricci flow,
see [T],[Lo],[B] and the reference therein. In [T] Topping
introduced $\mathcal{L}$-optimal transportation for Ricci flow. He
studied the behavior of Boltzmann-Shannon entropy along
$\mathcal{L}$-Wasserstein geodesic, and obtained natural monotonic
quantity from which the monotonicity of Perelman's
$\mathcal{W}$-entropy was recovered among other things. Lott [Lo]
showed the convexity of a certain entropy-like function using
Topping's work [T], as a result, he could reprove the monotonicity
of Perelman's reduced volume. In [B] Brendle proved a
Pr$\acute{e}$kopa-Leindler-type inequality for Ricci flow using [T],
from which he could also recover the monotonicity of Perelman's
reduced volume.

On the other hand, List [Li] considered an extended Ricci flow in
his thesis, and he generalized the monotonicity of Perelman's
$\mathcal{W}$-entropy to his flow. M$\ddot{u}$ller [M] studied
more general evolving closed manifolds $(M,g_{ij}(t))$ with the
metrics $g_{ij}(t)$ satisfying the equation

$\frac{\partial g_{ij}}{\partial t}=-2S_{ij}$,       $\  \    \   \
    \     \     \    \     \   $ (1.1)

where $\mathcal{S}=(S_{ij})$ is a symmetric tensor field of
(2,0)-type on $M$. He generalized the monotonicity of Perelman's
reduced volume to this flow satisfying a certain constraint
condition which will be stated later; see [M, Theorem 1.4].

In this note we will adapt Topping's  $\mathcal{L}$-optimal
transportation theory for Ricci flow to the  general flow (1.1). We
obtain some analogs of results of Topping, Lott and Brendle
mentioned above, and using this we can generalize the monotonicity
of List's (and hence also of Perelman's) $\mathcal{W}$-entropy, and
recover the monotonicity of M$\ddot{u}$ller's (and hence also of
Perelman's) reduced volume.

Now we consider the flow (1.1) backwards in time. Let $\tau$ be
some backward time parameter (i.e. $\tau=C-t$ for some constant $C
\in \mathbb{R}$). Consider the reverse flow

$\frac{\partial g_{ij}}{\partial \tau}=2S_{ij}(\tau)$, $ \    \
    \     \     \     \   \   $ (1.2)

 defined on a time interval including $[\tau_1, \tau_2]$ (with $0\leq \tau_1< \tau_2$).
 Following Perelman [P] and M$\ddot{u}$ller [M], we define the
 $\mathcal{L}$-length of a curve
 $\gamma:[\tau_1,\tau_2]\rightarrow M$ by

$\mathcal{L}(\gamma):=\int_{\tau_1}^{\tau_2} \sqrt{\tau}(S(\gamma
(\tau),\tau)+|\gamma'(\tau)|_{g(\tau)}^2)d\tau$,

where  $S$ is the trace of $\mathcal{S}$ ( w.r.t. $g(\tau)$). Then
we define the $\mathcal{L}$-distance by

$Q(x,\tau_1;y,\tau_2):=inf\{\mathcal{L}(\gamma)$ $|
\gamma:[\tau_1,\tau_2]\rightarrow M$ is smooth and $\gamma
(\tau_1)=x, \gamma(\tau_2)=y\}$.

Given two Borel probability measures $\nu_1,\nu_2$ viewed at times
$\tau_1$ and $\tau_2$ respectively , following [T] we define the
$\mathcal{L}$-Wasserstein distance by

$V(\nu_1,\tau_1;\nu_2,\tau_2):=inf \{\int_{M \times
M}Q(x,\tau_1;y,\tau_2)d\pi (x,y)$  $| \pi \in \Gamma
(\nu_1,\nu_2)\}$, $ \     \     \    \  $ (1.3)

where $\Gamma (\nu_1,\nu_2)$ is the space of Borel probability
measures on $M \times M$ with marginals $\nu_1$ and $\nu_2$.

 To state our theorems we
need to introduce a quantity in [M].
 Let $g(\tau)$ evolve by (1.2), and let $X \in \Gamma(TM)$ be a vector field on $M$. Set

$\mathcal{D}(\mathcal{S},X):=-\partial_{\tau}S-\triangle
S-2|S_{ij}|^2+4(\nabla_iS_{ij})X_j-2(\nabla_jS)X_j+2R_{ij}X_iX_j-2S_{ij}X_iX_j$.

 Our first result generalizes [T,
 Theorem 1.1] and a result of von Renesse  and Sturm [vRS]. As in [T] we refer to a family of smooth probability measures
$\nu(\tau)$ on $M$ as a diffusion if the density $u(\tau)$ relative
to the Riemannian volume measure $\mu (\tau)$ of $g(\tau)$ (i.e.
$d\nu (\tau)= u(\tau)d\mu (\tau)$) satisfies the equation

$\frac{\partial u}{\partial \tau}=\triangle u-Su$. $\    \   \
 \    \ $ (1.4)

\hspace *{0.4cm}

  {\bf Theorem 1.1} Given  $0 <\bar{\tau}_1< \bar{\tau}_2$, suppose
  that
  $(M,g(\tau))$ is a closed, $n$-dimensional  manifold  evolving by (1.2), for $\tau$ in some
 open interval containing $[\bar{\tau}_1,\bar{\tau}_2]$, such that the
 quantity $\mathcal{D}(\mathcal{S},X)$ is nonnegative for all
 vector fields $X \in \Gamma (TM)$ and all times for which the
 flow exists. Let $\nu_1(\tau)$ and $\nu_2(\tau)$ be two
 diffusions (as defined above) for $\tau$ in some neighbourhoods of $\bar{\tau}_1$ and $
  \bar{\tau}_2$ respectively. Set
  $\tau_1=\tau_1(s):=\bar{\tau}_1e^s,
  \tau_2=\tau_2(s):=\bar{\tau}_2e^s$, and define the renormalized
  $\mathcal{L}$- Wasserstein distance by

  $\Theta(s):=2(\sqrt{\tau_2}-\sqrt{\tau_1})V(\nu_1(\tau_1),\tau_1;\nu_2(\tau_2),\tau_2)
  -2n(\sqrt{\tau_2}-\sqrt{\tau_1})^2$

  for $s$ in a neighbourhood of 0 such that $\nu_i(\tau_i(s))$ are
  defined ($i=1,2$).

  Then $\Theta(s)$ is a weakly decreasing function of $s$.

\hspace *{0.4cm}

The constraint condition on  $\mathcal{D}(\mathcal{S},X)$ in Theorem
1.1 is the same as that appeared in [M, Theorem 1.4] mentioned
above. As pointed out in [M], it is satisfied, for example, by the
static manifolds with nonnegative Ricci curvature, by Hamilton's
Ricci flow, by List's flow ([Li]), by the Ricci flow coupled with
harmonic map heat flow  introduced by M$\ddot{u}$ller in his thesis
( cf. [M]), and by mean curvature flow in an ambient Lorentzian
manifold with nonnegative sectional curvature.

\hspace *{0.4cm}

Our second result generalizes [Lo, Theorem 1].

\hspace *{0.4cm}

 {\bf Theorem 1.2}  Given  $0 <\tau_1< \tau_2$, suppose
  that
  $(M,g(\tau))$ is a connected closed manifold evolving by (1.2), for $\tau $ in some open interval including $[\tau_1,\tau_2]$, such that the
 quantity $\mathcal{D}(\mathcal{S},X)$ is nonnegative for all
 vector fields $X \in \Gamma (TM)$ and all times for which the
 flow exists. Let
 $\mathcal{V}_{\tau}(\tau \in[\tau_1, \tau_2])$ be an $\mathcal{L}$-Wasserstein
 geodesic,induced by a potential $\varphi:M\rightarrow \mathbb{R}$, with $\mathcal{V}_{\tau_1}$
 and $\mathcal{V}_{\tau_2}$
 both absolutely continuous probability measures. Set $\phi(y,\tau):=\frac{1}{2\sqrt{\tau}}$ inf$_{x \in M}[Q(x,\tau_1;y,\tau)-\varphi(x)]$ for $y \in M$
 and $\tau \in [\tau_1,\tau_2]$. Then $E(\mathcal{V}_{\tau})+
 \int_M \phi(\cdot,\tau)d\mathcal{V}_{\tau}+\frac{n}{2}$ln$\tau$ is convex in the variable $\tau^{-1/2}$.

\hspace *{0.4cm}

For the definition of $\mathcal{L}$-Wasserstein
 geodesic see the paragraph following Theorem 2.14 in [T], cf. also the paragraph following our Theorem 2.1.
 Also note that here $E(\mathcal{V}_{\tau})$ is the Boltzmann-Shannon entropy of
 $\mathcal{V}_{\tau}$ (cf. Section 2).

\hspace *{0.4cm}

Our third theorem generalizes [B, Theorem 2] and a result in [CMS].
Note that we do not assume that $M$ is compact in this theorem.

\hspace *{0.4cm}

{\bf Theorem 1.3}  Given  $0 <\tau_1< \tau_2$, suppose
  that
  $(M,g(\tau))$ is a  complete  manifold evolving by (1.2), for $\tau $ in some open interval including $[\tau_1,\tau_2]$,
   with the sectional curvature and $S_{ij}$ uniformly bounded
  in compact time intervals, and such that the
 quantity $\mathcal{D}(\mathcal{S},X)$ is nonnegative for all
 vector fields $X \in \Gamma (TM)$ and all times for which the
 flow exists. Fix $\bar{\tau} \in (\tau_1,\tau_2)$, and write

 $\frac{1}{\sqrt{\bar{\tau}}}=\frac{1-\lambda}{\sqrt{\tau_1}}+\frac{\lambda}{\sqrt{\tau_2}}$,

 for some $0< \lambda <1$. Let $ u_1$, $u_2$, $v: M \rightarrow  \mathbb{R}$
 be nonnegative measurable functions such that

 $(\frac{\bar{\tau}}{\tau_1^{1-\lambda}
 \tau_2^{\lambda}})^{\frac{n}{2}}v(\gamma(\bar{\tau}))\geq$
 exp$(-\frac{1-\lambda}{2\sqrt{\tau_1}}Q(\gamma
 (\tau_1),\tau_1;\gamma (\bar{\tau}),\bar{\tau}))$ $ u_1(\gamma
 (\tau_1))^{1-\lambda}$

 $\cdot$ exp$(\frac{\lambda}{2\sqrt{\tau_2}}Q(\gamma
 (\bar{\tau}),\bar{\tau};\gamma (\tau_2),\tau_2))$ $ u_2(\gamma
 (\tau_2))^{\lambda}$

 for each minimizing $\mathcal{L}$-geodesic $\gamma :
 [\tau_1,\tau_2] \rightarrow M$. Then

 $\int_Mvd\mu(\bar{\tau})\geq
 (\int_Mu_1d\mu(\tau_1))^{1-\lambda}(\int_Mu_2d\mu(\tau_2))^{\lambda}$.

\hspace *{0.4cm}

In Section 2 we give the proof of our theorems which relies heavily
on Topping [T]. In Section 3 we give some applications of our
theorems (following Topping and Brendle).

\section {Proof of Theorems}

Part of Topping's $\mathcal{L}$-optimal transportation theory for
Ricci flow [T] extends to the  general flow (1.1) without any
change. In particular, virtually all theorems in [T, Section 2] hold
in our more general situation. We just state the following

\hspace *{0.4cm}

{\bf Theorem 2.1}(cf. [T, Section 2, in particular Theorem 2.14])
Given $0 <\tau_1< \tau_2$, suppose
  that
  $(M,g(\tau))$ is a  closed manifold evolving by (1.2), for $\tau $ in some open interval including $[\tau_1,\tau_2]$.
   Suppose that $\nu_1$ and
$\nu_2$ are absolutely continuous probability measures (w.r.t. (any)
volume measure). Then there exists an optimal transference plan
$\pi$ in (1.3) which is given by the push-forward of $\nu_1$ under
the map $x\mapsto (x,F(x))$, where $F: M\rightarrow M$ is a Borel
map defined by

$F(x):=\mathcal{L}_{\tau_1,\tau_2}exp_x(-\frac{\nabla
\varphi(x)}{2})$,           $\   \   \    \   \   \  $  (2.1)

at points of differentiability of some reflexive function $\varphi:
M\rightarrow R$, where the gradient is w.r.t. $g(\tau_1)$.

Moreover, there exists a Borel set $K\subset M$ with $\nu_1(K)=1$,
such that for each $x \in K$, $\varphi$ admits a Hessian at $x$, and

$f_{\tau_1}(x)=f_{\tau_2}(F(x))$ det $(dF)_x \neq 0$,   $\   \   \ \
   \  $  (2.2)

 where $f_{\tau_i}$ is the densities defined by
$d\nu_i=f_{\tau_i}d\mu(\tau_i)$ for $i=1,2$.

\hspace *{0.4cm}

 As in [T], we refer to $\mathcal{V}_{\tau}:=(F_{\tau})_\sharp(\nu_1)$
as an $\mathcal{L}$-Wasserstein geodesic, where
$F_{\tau}:M\rightarrow M$ is a Borel map defined by

$F_{\tau}(x):=\mathcal{L}_{\tau_1,\tau}exp_x(-\frac{\nabla
\varphi(x)}{2})$

at points of differentiability of $\varphi$ (as in the above
theorem) for $\tau \in [\tau_1, \tau_2]$.

\hspace *{0.4cm}

{\bf Remark 2.2} Theorem 2.1 extends to noncompact case with
suitable  modifications. More precisely, when $M$ is noncompact, one
imposes in addition the conditions that $S_{ij}$ is uniformly
bounded (in compact time intervals)  and that
$V(\nu_1,\tau_1;\nu_2,\tau_2)$ is finite, then the results in
Theorem 2.1 still hold with the gradient in (2.1) and the
differential in (2.2) replaced by an approximate gradient and an
approximate differential respectively, and Hessian replaced by
approximate Hessian.( Of course, $\varphi$ need not be reflexive any
more.) For more details, one can consult [FF], [F] and [V].
Moreover, in noncompact case, even if one does not impose the
finiteness condition on $V(\nu_1,\tau_1;\nu_2,\tau_2)$, one can
still say something, cf. [F] and [V].

\hspace *{0.4cm}

Note that M$\ddot{u}$ller [M] has established some properties of
$\mathcal{L}$-geodesics and $L$-function in our situation.

As in [M], we introduce

$\mathcal{H}(\mathcal{S},X):=-\partial_{\tau}S-\frac{1}{\tau}S-2X(S)+2\mathcal{S}(X,X)$.

The following lemma generalizes [T, Lemma 3.1].

 \hspace *{0.4cm}

{\bf Lemma 2.3} Let $\gamma:[\tau_1,\tau_2]\rightarrow M$ be an
$\mathcal{L}$-geodesic, and $\{Y_i(\tau)\}_{i=1,...,n}$ be a set of
$\mathcal{L}$-Jacobi fields along $\gamma$ which form a basis of
$T_{\gamma(\tau)}M$ for each $\tau \in [\tau_1,\tau_2]$, with
$\{Y_i(\tau_1)\}$ orthonormal and $\langle D_\tau Y_i, Y_j\rangle$
symmetric in $i$ and $j$ at $\tau=\tau_1$. Define
$\alpha:[\tau_1,\tau_2] \rightarrow \mathbb{R}$ by
 $\alpha
(\tau)=-\frac{1}{2}$ln det$\langle Y_i(\tau),
Y_j(\tau)\rangle_{g(\tau)}$, and write $\sigma=\sqrt{\tau}$, then we
have

$\frac{d^2\alpha}{d\sigma^2}=4\sqrt{\tau}\frac{d}{d\tau}(\sqrt{\tau}\frac{d\alpha}{d\tau})\geq
2\tau (\mathcal{H}(\mathcal{S},X)+ \mathcal{D}(\mathcal{S},X))$,

and

$\frac{d^2(\sigma\alpha)}{d\sigma^2}=4\frac{d}{d\tau}(\tau^{\frac{3}{2}}\frac{d\alpha}{d\tau})\geq
2\tau^{\frac{3}{2}} (\mathcal{H}(\mathcal{S},X)+
\mathcal{D}(\mathcal{S},X))-n\tau^{-\frac{1}{2}}$,

where $X=\gamma'(\tau)$.

\hspace *{0.4cm}

{\bf Proof } The proof follows closely that of Topping [T, Lemma
3.1] with some necessary modifications. From the
$\mathcal{L}$-geodesic equation in [M] we can derive  the
$\mathcal{L}$-Jacobi equation for $Y(\tau)$

$D_\tau^2 Y:=D_\tau(D_\tau Y)$= $-R(X,Y)X+\frac{1}{2}\nabla_Y(\nabla
S)-\nabla_Y\widetilde{\mathcal{S}}(X)-2\widetilde{\mathcal{S}}(D_\tau
Y)-\frac{1}{2\tau}D_\tau Y+
 \nabla_X \widetilde{\mathcal{S}}(Y)-[\nabla
\mathcal{S}(\cdot,X,Y)]^\sharp$.

Here $\widetilde{\mathcal{S}}$ is $\mathcal{S}$ viewed as an
endomorphism (i.e. a (1,1)- tensor), other conventions are from
[T].

 Consider the solution
$e_i\in \Gamma(\gamma^*(TM)), i=1,...,n,$ of the ODE

$D_\tau e_i+\widetilde{\mathcal{S}}(e_i)=0,$

with initial condition $e_i(\tau_1):=Y_i(\tau_1)$. Write
$Y_j(\tau)= A_{kj}e_k(\tau)$ for a $\tau$-dependent $n \times n$
matrix $A$. Then we have

$A_{ij}'=\langle D_\tau Y_j,e_i\rangle+
A_{kj}\mathcal{S}(e_k,e_i)$,

and

$A_{ij}''=\langle D_\tau^2Y_j,e_i\rangle +2
A_{kj}'\mathcal{S}(e_k,e_i)+ A_{kj}\langle D_\tau
(\widetilde{\mathcal{S}}(e_k)),e_i\rangle$.

Using the $\mathcal{L}$-Jacobi equation we get that

$\langle D_\tau^2Y_j,e_i\rangle=
A_{kj}[-Rm(X,e_k,X,e_i)+\frac{1}{2}$ Hess$(S)(e_i,e_k)+ \nabla_X
\mathcal{S}(e_i,e_k)-\langle
\nabla_{e_k}\widetilde{\mathcal{S}}(X),e_i\rangle- \langle
\nabla_{e_i}\widetilde{\mathcal{S}}(X),e_k\rangle + 2 \langle
\widetilde{\mathcal{S}}^2(e_k),e_i\rangle+\frac{1}{2\tau}\mathcal{S}(e_i,e_k)]
-2 A_{kj}'\mathcal{S}(e_k,e_i)-\frac{1}{2\tau}A_{ij}'.$

We  also have

$\langle D_\tau
(\widetilde{\mathcal{S}}(e_k)),e_i\rangle=\frac{\partial
\mathcal{S}}{\partial \tau}(e_i,e_k)+\nabla_X
\mathcal{S}(e_i,e_k)-3\langle
\widetilde{\mathcal{S}}^2(e_k),e_i\rangle$.

Then we get that

$A''+\frac{1}{2\tau}A'=MA$,

where $M$ is the  $\tau$-dependent $n \times n$ symmetric matrix
given by

$M_{ik}=-Rm(X,e_k,X,e_i)+\frac{1}{2}$ Hess$(S)(e_i,e_k)+ 2\nabla_X
\mathcal{S}(e_i,e_k)-\langle
\nabla_{e_k}\widetilde{\mathcal{S}}(X),e_i\rangle- \langle
\nabla_{e_i}\widetilde{\mathcal{S}}(X),e_k\rangle - \langle
\widetilde{\mathcal{S}}^2(e_k),e_i\rangle+\frac{1}{2\tau}\mathcal{S}(e_i,e_k)+\frac{\partial
\mathcal{S}}{\partial \tau}(e_i,e_k).$

Using [M, Lemma 1.6], we see that the trace of $M$ is

tr$M=-\frac{1}{2}(\mathcal{H}(\mathcal{S},X)+
\mathcal{D}(\mathcal{S},X))$.

Now  define $B:=\frac{dA}{d\tau}A^{-1}$, then similarly as in [T],
we have

$\tau^{-1/2}\frac{d}{d\tau}(\sqrt{\tau}\frac{d\alpha}{d\tau})=$ tr
$B^2+\frac{1}{2}(\mathcal{H}(\mathcal{S},X)+
\mathcal{D}(\mathcal{S},X))$,

and

$\tau^{-\frac{3}{2}}\frac{d}{d\tau}(\tau^{\frac{3}{2}}\frac{d\alpha}{d\tau})=$
tr $(B-\frac{1}{2\tau}I)^2+\frac{1}{2}(\mathcal{H}(\mathcal{S},X)+
\mathcal{D}(\mathcal{S},X))-\frac{n}{4\tau^2}.$

Similarly as  in [T], one can show that $B$ is symmetric, and our
result follows.

\hspace *{0.4cm}

Now we begin to study the behavior of Boltzmann-Shannon entropy
along a $\mathcal{L}$-Wasserstein geodesic. Recall that the
Boltzmann-Shannon entropy of a probability measure $fd\mu$ is
defined by

$E(fd\mu)=\int_Mf$ln$fd\mu$,

where $\mu$ is Riemannian volume measure, and $f$ is a
 reasonably regular weakly positive function on $M$. As before we set $\sigma=\sqrt{\tau}$. Then we have
the following lemma which generalizes [T,Lemma 3.2].

\hspace *{0.4cm}

{\bf Lemma 2.4} Let $(M,g(\tau))$ be as in Theorem 2.1. Let
 $\mathcal{V}_{\tau}(\tau \in[\tau_1, \tau_2])$ be an $\mathcal{L}$-Wasserstein
 geodesic, induced by a potential $\varphi:M\rightarrow \mathbb{R}$, with $\mathcal{V}_{\tau_1}$
  and $\mathcal{V}_{\tau_2}$
 both absolutely continuous probability measures, and write $d\mathcal{V}_{\tau}=f_{\tau}d\mu(\tau)$ where $\mu(\tau)$ is the volume
 measure of $g(\tau)$. Then for all $\tau \in[\tau_1, \tau_2]$, we have $f_\tau \in L$ln$L(\mu(\tau))$,
 and the function $
 E(\mathcal{V}_{\tau})$ is semiconvex  in $\tau$ and satisfies, for almost all
 $\tau \in[\tau_1, \tau_2]$ (where $\sigma \mapsto E(\mathcal{V}_{\tau})$
 admits a second derivative in the sense of Alexandrov)

 $\frac{d^2}{d\sigma^2}E(\mathcal{V}_{\tau})=
 4\sqrt{\tau}\frac{d}{d\tau}(\sqrt{\tau}\frac{dE(\mathcal{V}_{\tau}) }{d\tau})\geq
2\tau \int_M(\mathcal{H}(\mathcal{S},X(\tau))+
\mathcal{D}(\mathcal{S},X(\tau)))d\mathcal{V}_{\tau_1}$,

and

$\frac{d^2}{d\sigma^2}(\sigma
E(\mathcal{V}_{\tau}))=4\frac{d}{d\tau}(\tau^{\frac{3}{2}}\frac{dE(\mathcal{V}_{\tau})
}{d\tau})\geq 2\tau^{\frac{3}{2}}
\int_M(\mathcal{H}(\mathcal{S},X(\tau))+
\mathcal{D}(\mathcal{S},X(\tau)))d\mathcal{V}_{\tau_1}-n\tau^{-\frac{1}{2}}$,

where $X(\tau)$, at a point $x \in M$ where $\varphi$ admits a
Hessian, is $\gamma'(\tau)$, for $\gamma:[\tau_1,
\tau_2]\rightarrow M$ the minimizing $\mathcal{L}$-geodesic from
$x$ to $F(x)$. Moreover, the one-sided derivatives of
$E(\mathcal{V}_{\tau})$ at $\tau_1$ and $\tau_2$ exist, with

$\frac{d}{d\tau}|_{\tau_1}E(\mathcal{V}_{\tau})\geq
-\int_M(S(\cdot,\tau_1)+\langle \frac{\nabla
\varphi}{2\sqrt{\tau_1}},\nabla
$ln$f_{\tau_1}\rangle_{g(\tau_1)})d\mathcal{V}_{\tau_1}$.

\hspace *{0.4cm}

{\bf Proof} Using Theorem 2.1 and Lemma 2.3 one can proceed exactly
as in [T].

\hspace *{0.4cm}

Now suppose $g(\tau)$ is defined on
$(\hat{\tau}_1,\hat{\tau}_2)\supset [\tau_1,\tau_2]$, where
$\hat{\tau}_1 > 0$. As in  [T], let
$\Upsilon:=\{(x,\tau_a;y,\tau_b)|x,y \in M$ and $\hat{\tau}_1<
\tau_a < \tau_b < \hat{\tau}_2\}$. Suppose $(x,\tau_1;y,\tau_2)\in
\Upsilon \setminus \mathcal{L}Cut$, let
$\gamma:[\tau_1,\tau_2]\rightarrow M$ be the minimizing
$\mathcal{L}$-geodesic from $x$ to $y$, and write
$X(\tau)=\gamma'(\tau)$ as before. Following [P],[T] and [M], define

$\mathcal{K}=\mathcal{K}(x,\tau_1,y,\tau_2):=\int_{\tau_1}^{\tau_2}\tau^{\frac{3}{2}}
\mathcal{H}(\mathcal{S},X(\tau))d\tau$.

Then we have the following result which generalizes [T. Corollary
3.3].

\hspace *{0.4cm}

{\bf Corollary 2.5} Let the hypothesis of Lemma 2.4 still hold, and
assume further that the
 quantity $\mathcal{D}(\mathcal{S},X)$ is nonnegative for all
 vector fields $X \in \Gamma (TM)$ and all times for which the
 flow exists. Then

 $\int_{M\times M}(\mathcal{K}-2\tau_1^{3/2}S(x,\tau_1)-\tau_1\langle
 \nabla_1
Q,\nabla $
ln$f_{\tau_1}(x)\rangle_{g(\tau_1)}+2\tau_2^{3/2}S(y,\tau_2)$

$-\tau_2\langle
 \nabla_2
Q,\nabla $ln$f_{\tau_2}(y)\rangle_{g(\tau_2)})d\pi(x,y)$

$\leq n(\sqrt{\tau_2}-\sqrt{\tau_1}),$

where $\nabla_1Q$ denotes the gradient of $Q$ w.r.t. its $x$
argument and w.r.t. $g(\tau_1)$, $\nabla_2Q$ denotes the gradient of
$Q$ w.r.t. its $y$ argument and w.r.t. $g(\tau_2)$,and $\pi$ is the
optimal transference plan from $\mathcal{V}_{\tau_1}$ to
$\mathcal{V}_{\tau_2}$ (for $\mathcal{L}$-optimal transportation).

\hspace *{0.4cm}

The following result generalizes [T, Lemma A.6].

\hspace *{0.4cm}

{\bf Lemma 2.6} Under the flow (1.2), we have

 $\tau_2\frac{\partial Q}{\partial \tau_2}+\tau_1\frac{\partial Q}{\partial
 \tau_1}=2\tau_2^{\frac{3}{2}}S(y,\tau_2)-2\tau_1^{\frac{3}{2}}S(x,\tau_1)+\mathcal{K}-\frac{1}{2}Q$.

\hspace *{0.4cm}

{\bf Proof} Similarly as [T, (A.4) and (A.5)],we have

$\frac{\partial Q}{\partial
\tau_1}(x,\tau_1;y,\tau_2)=\sqrt{\tau_1}(|X(\tau_1)|^2-S(x,\tau_1));$
$\nabla_1Q(x,\tau_1;y,\tau_2)=-2\sqrt{\tau_1}X(\tau_1)$,

and

$\frac{\partial Q}{\partial
\tau_2}(x,\tau_1;y,\tau_2)=\sqrt{\tau_2}(S(y,\tau_2))-|X(\tau_2)|^2);$
$\nabla_2Q(x,\tau_1;y,\tau_2)=2\sqrt{\tau_2}X(\tau_2)$.

Similarly as [T, (A.9)], we have

$\tau_2^{\frac{3}{2}}(S(y,\tau_2)+|X(\tau_2)|^2)-\tau_1^{\frac{3}{2}}(S(x,\tau_1)+|X(\tau_1)|^2)
=-\mathcal{K}(x,\tau_1,y,\tau_2)+\frac{1}{2}Q(x,\tau_1;y,\tau_2)$.

(cf. also [M].)

Then the lemma follows.

\hspace *{0.4cm}

Finally, Theorem 1.1 follows from Corollary 2.5 and Lemma 2.6 (cf.
[T, Section 4]).

\hspace *{0.4cm}

For the proof of Theorem 1.2, we follow closely [Lo].

From [T, Lemma 2.4] we can derive

$\tau^{\frac{3}{2}}\frac{d}{d\tau}\phi(\gamma(\tau))=-\frac{1}{2}\sqrt{\tau}\phi(\gamma(\tau))+
\frac{1}{2}\tau^{\frac{3}{2}}(S(\gamma(\tau),\tau)+|X(\tau)|^2)$.

From [M] we have

$\frac{d}{d\tau}(S(\gamma(\tau),\tau)+|X(\tau)|^2)=-\mathcal{H}(\mathcal{S},X)-
\frac{1}{\tau}(S(\gamma(\tau),\tau)+|X(\tau)|^2)$.

It follows that

$(\tau^{\frac{3}{2}}\frac{d}{d\tau})^2\phi(\gamma(\tau))=-\frac{1}{2}\tau^3\mathcal{H}(\mathcal{S},X)$.

Combining with the condition
$\mathcal{V}_{\tau}=(F_{\tau})_\sharp(\nu_1)$, the equation above
implies

$(\tau^{\frac{3}{2}}\frac{d}{d\tau})^2\int_M\phi(\tau)d\mathcal{V}_{\tau}=
-\frac{1}{2}\tau^3 \int_M\mathcal{H}(\mathcal{S},\nabla
\phi(\tau))d\mathcal{V}_{\tau}$.

Combining the equation above  with Lemma 2.4 and the assumption on
$\mathcal{D}(\mathcal{S},X)$, we get Theoren 1.2.

\hspace *{0.4cm}

 To prove Theorem 1.3, it suffices to prove the case
that $u_1, u_2$ have compact support, since then the general case
will follows by an approximate technique as in [CMS]. Now one
proceeds as in [B]. A key step is to prove that under our assumption
on $\mathcal{D}(\mathcal{S},X)$, one has

 $\tau^{-\frac{3}{2}}\frac{d}{d\tau}[\tau^{\frac{3}{2}}\frac{d}{d\tau}(\frac{n}{2}$
 ln $\tau +\frac{1}{2}\tau^{-\frac{1}{2}}Q(x,\tau_1;F_{\tau}(x),\tau)-$ ln det $(dF_{\tau})_x)] \geq 0$

 as in [B]. This can be proved by using Lemma 2.3, similarly as in the proof of Theorem 1.2.

\section {Some applications}

For closed manifold $(M,g(\tau))$ evolving by (1.2) and a solution
$u$ of (1.4) we introduced the $\mathcal{W}$-entropy as in
[P],[Li],

$\mathcal{W}:=\int_M[\tau(S+|\nabla f|^2)+f-n](4\pi
\tau)^{-n/2}e^{-f}d\mu$,

where $f$ is defined by $u=(4\pi \tau)^{-n/2}e^{-f}$.

Then we have the following

\hspace *{0.4cm}

 {\bf Theorem 3.1} Assume that the quantity $\mathcal{D}(\mathcal{S},X)$ is nonnegative for all
 vector fields $X \in \Gamma (TM)$ and all times for which the
 flow exists. Then $\frac{d\mathcal{W}}{d\tau}\leq 0$.

\hspace *{0.4cm}

{\bf Proof} Theorem 3.1 follows easily from Theorem 1.1 and a result
which generalizes [T, Lemma 1.3] (with $S$ replacing $R$ in
[T,(1.7)]) and whose proof is a minor modification of that of [T,
Lemma 1.3].

\hspace *{0.4cm}

{\bf Remark 3.2} Some special cases of Theorem 3.1 appeared in [P],
[N] and [Li]. Of course, one can also prove Theorem 3.1 by a direct
computation as in these references.

\hspace *{0.4cm}

As in [T, Section 1.3], Theorem 1.1 also implies the monotonicity of
the enlarged length which generalizes the corresponding result of
Perelman [P]. More precisely, as in [P], consider
$L(y,\tau):=Q(x,0;y,\tau)$  for fixed $x \in M$ and
$\bar{L}(y,\tau):=2\sqrt{\tau}L(y,\tau)$. Then we have the following

\hspace *{0.4cm}

{\bf Theorem 3.3} Assume that the quantity
$\mathcal{D}(\mathcal{S},X)$ is nonnegative for all
 vector fields $X \in \Gamma (TM)$ and all times for which the
 flow exists. Then the minimum over $M$ of
$\bar{L}(\cdot,\tau)-2n\tau$ is a weakly decreasing function of
$\tau$.

\hspace *{0.4cm}

The following theorem extends a theorem in [P], and also extends a
theorem in [M] to the noncompact case.

\hspace *{0.4cm}

 {\bf Theorem  3.4} Suppose
  that
  $(M,g(\tau))$ is a  complete  manifold evolving by (1.2),  with the sectional curvature and $S_{ij}$ uniformly bounded
  in compact time intervals, and such that the
 quantity $\mathcal{D}(\mathcal{S},X)$ is nonnegative for all
 vector fields $X \in \Gamma (TM)$ and all times for which the
 flow exists. Then the  reduced volume (as defined in [P], [M]) is
 nonincreasing in $\tau$.

\hspace *{0.4cm}

{\bf Proof}.  This is a corollary of Theorem 1.3, cf. [B, Section
3].

\hspace *{0.4cm}

 {\bf Remark 3.5} Our $\mathcal{L}$-length is the
same as $\mathcal{L}_b$-length in [M], and corresponds to
$\mathcal{L}_-$-length in [Lo]. One can also develop a parallel
theory of $\mathcal{L}_f$-(or $\mathcal{L}_+$-) and
$\mathcal{L}_0$-optimal transportation respectively as in [Lo].

\bibliographystyle{amsplain}

\hspace *{0.4cm}

{\bf Reference}

\bibliography{1}[B] S. Brendle, A Pr$\acute{e}$kopa-Leindler-type inequality for Ricci flow, arXiv:0907.3726.

\bibliography{2}[CMS] D. Cordero-Erausquin, R.J. McCann and M.
Schmuckenschl$\ddot{a}$ger, A Riemannian interpolation inequality
$\grave{a}$ la Borell, Brascamp and Lieb, Invent. Math. 146 (2001),
219-257.

\bibliography{3}[FF] A. Fathi and A. Figalli, Optimal transportation
on non-compact manifolds, arXiv:0711.4519, to appear in Israel J.
Math.

\bibliography{4}[F]  A. Figalli, Existence, uniqueness and
regularity of optimal transport maps, SIAM J. Math. Anal. 39 (2007),
126-137.

\bibliography{5}[Li] B. List, Evolution of an extended Ricci flow system, Comm. Anal. Geom. 16 (2008), 1007-1048.

\bibliography{6}[Lo] J. Lott, Optimal transport and Perelman's reduced volume, arXiv:0804.0343v2, to appear in Cal. Var. PDE.

\bibliography{7}[M] R. M$\ddot{u}$ller, Monotone volume formulas for geometric flows, arXiv:0905.2328,
to appear in J. Reine Angew. Math.

\bibliography{8}[N] L. Ni, The entropy formula for linear heat
equation, J. Geom. Anal. 14 (2004), 85-98; Addenda, 14
(2004),369-374.

\bibliography{9}[P] G. Perelman, The entropy formula for the Ricci flow and its
geometric applications, arXiv:math.DG/0211159.

\bibliography{10}[vRS] M.-K. von Renesse  and K.-T. Sturm, Transport
inequalities, gradient estimates, entropy and Ricci curvature, Comm.
Pur Appl. Math. 58 (2005), 923-940.

\bibliography{11}[T] P. Topping, $\mathcal{L}$-optimal transportation
for Ricci flow, available online at Topping's home page, and to
appear in J. Reine Angew. Math.

\bibliography{12}[V] C. Villani, Optimal transport, old and
new, Springer-Verlag, 2009.

\end{document}